\documentclass[12pt, a4paper]{article}
\usepackage{amsmath,amssymb,amscd,amsthm}

\newcommand{\Herpos}{\dot{P}_m}

\newcommand{\QFmet}[2]{(\!( {#1} , {#2} )\!)^{QF}}
\newcommand{\SPHmet}[2]{(\!( {#1} , {#2} )\!)^{Sph}}
\newcommand{\Dmet}[2]{(\!( {#1} , {#2} )\!)^{D}}
\newcommand{\Smet}[2]{(\!( {#1} , {#2} )\!)^{S}}

\newcommand{\tr}{\mathrm{tr}}
\newtheorem{definition}{Definition}[section]
\newtheorem{theorem}[definition]{Theorem}
\newtheorem{lemma}[definition]{Lemma}

\newtheorem{remark}{Remark}

\begin{document}
\begin{center}
{\Large
A gradient system on the quantum information space
}
\\ \bigskip
{\Large
realizing the averaged learning equation of Hebb type
}
\\ \bigskip
{\Large
for the principal component analyzer
}
\end{center}
\bigskip
\begin{center}
{\large
Yoshio Uwano
and
Hiromi Yuya
}
\\ \bigskip \bigskip
Department of Complex Systems
\\
Future University Hakodate
\\
Kameda Nakano-cho 116-2, Hakodate
\\
041-8655, Japan
\\ \medskip
e-mail: uwano@fun.ac.jp, g2108042@fun.ac.jp
\end{center}
\smallskip
\begin{abstract}
The averaged learning equation (ALEH) applicable to the principal component
analyzer is studied from both quantum information geometry and
dynamical system viewpoints. On the quantum information space (QIS),
the space of regular density matrices endowed with the quantum
SLD-Fisher metric, a gradient system is given as an extension of
the ALEH; on the submanifold, consisting of the diagonal matrices, of
the QIS, the gradient flow coincides with the ALEH up to a local
diffeomorphism.

\end{abstract}
\section{Introduction}
Quantum computing has been widely known to be a research area
in rapid-rate progress. It has been recognized also as
an interdisciplinary research area in which numbers of researchers
with a variety of backgrounds are working.
An ambition of breaking theoretical boundaries of binary computing
has been a great driving force behind many investigations
to discover quantum algorithms beyond the boundaries.
As celebrated examples beyond the boundaries,
Shor's algorithms and Grover's one are well-known for
the discrete logarithms, the prime factorization and
the data search, respectively \cite{Sh,G}.
To trace the history of quantum computing, see \cite{NC} for example. 
\par
On turning to algorithms in convensional sense of computing,
there exist various excellent algorithms developed in engineering
and systems science. It might happen that some of them admit a similar
mathematical structure: In 1990's, Nakamura revealed integrability,
Lax-type structure and gradient-system structure in
a matrix-eigenvalue computing \cite{Nev},
the gradient system on the space of multinomial distributions
\cite{Nstat}, the Karmarkar flow for linear programming \cite{NKar}
and an averaged learning equation of Hebb type \cite{NHebb}.
\par
In a series of papers by the authors \cite{U,UHI,UY,UY2},
a counterpart to Nakamura's gradient-system structure was succesfully
found in the quantum information space (QIS) for 
the gradient system on the space of multinomial distributions
\cite{Nstat} and the Karmarkar flow for linear programming \cite{NKar}.
The QIS in this paper stands for the space of regular density matrices
endowed with the quantum SLD-Fisher metric.
The former system is realized as the gradient system on the QIS associated
with the negative von Neumann entropy, the latter with a trace of the square
of density matrix times a cost-coefficient matrix.
As a continuation of those papers of the authors,
the aim of the present paper is to construct a counterpart of the
averaged learning equation of Hebb type for the principal component analyzer
\cite{NHebb,Oja}, which takes the gradient-system form in the QIS.
The present paper is placed in an interdisciplinary research area
of dynamical systems on the QIS, algorithms and applied differential geometry.
The result would be expected to be a clue to realize algorithms in the QIS
in mathematical sense first, and in a more realistic physical sense in future.
In what follows, the contents of the present paper are outlined.
\par
Section 2 is the preliminaries for the averaged learning equation
of Hebb type (ALEH). The ALEH is derived from a synaptic neuron model
and its gradient-system form is introduced.
Section 3 is devoted to geometric devices for realizing the ALEH
on the QIS. In order to transfer the ALEH on a dense submanifold
${\cal S}_m$ of the sphere to that on the submanifold ${\cal D}_m$
of the QIS consisting of diagonal matices, an immersion of ${\cal S}_m$
to ${\cal D}_m$ is introduced. On regarding the immersion as a local
diffeomorphism, the ALEH is understood to be transferred
to a multi-fold copy of the ALHE. It is worth
pointing out that the geometric devices presented in the present paper
are quite different from those in the papers \cite{UY,UY2} for the
Karmarkar flow. Section~4 is the core part of the present
paper, where the gradient system on the QIS realizing the ALEH is given
explicitly. Section~5 is for concluding remarks.
\section{Preliminary: The ALEH}
In this section, we review an averaged learning equation
of Hebb-type (ALEH) following Nakamura \cite{NHebb} and Oja \cite{Oja}. 
The gradient-system form derived by Nakamura \cite{NHebb} of the ALEH
is also reviewed on the $(m-1)$-dimensional unit sphere.
\subsection{Learning in the synaptic neuron model}
Let us introduce a vector-valued variable
$X=( X_{1}, X_{2}, \cdots, X_{m})^{T} \in {\bf R}^m$ to
express the $m$ presynaptic signals and a scalar-valued variable
$Y \in {\bf R}$
to express the postsynaptic signal; the values at the time $s$
of $X$ and $Y$ are described as $X(s)$ and $Y(s)$, respectively.
The synaptic neuron model dealt with in this paper starts with the
following timewise linear relation
\begin{eqnarray}
\label{output-signal}
Y(s) =W^{T}(s) X(s) 
\end{eqnarray}
between $X$ and $Y$, where
$W(s)=(W_{1}(s), W_{2}(s), \dots, W_{m}(s))^{T} \in {\bf R}^m$
stands for the vector of coupling strengths of the neuron.
Throughtout this paper, the superscript ${}^T$ indicates the
transpose operation.
\par
According to Hebb's hypothesis \cite{Oja,Hebb}, 
learning in the synaptic neuron models amounts to updating efficacies
of extracting inputs with high probability:
In the model (\ref{output-signal}), the vector of coupling strengths
$W(s)$ is understood to be updated through a recurrence relation along
with repetitive inputs of signals into neuron.
Following Oja \cite{Oja}, we consider the discrete-time recurrence relation
\begin{eqnarray}
\label{o-rule1}
W(s+1) = \frac{W(s) + \eta Y(s) X(s) }{\| W(s) + \eta Y(s) X(s)\|}
\end{eqnarray}
for the coupling strength $W(s)$ in (\ref{output-signal}),
where $\eta$ is a positive constant indicating the learning rate
and the symbol $\| \cdot \|$ stands for the standard
Euclidean norm of vectors in ${\bf R}^m$.
As a significant characteristic of
the recurrence relation (\ref{o-rule1}), the following
norm-preserving property is worth pointed out; 
\begin{equation}
\label{norm-preserve}
\| W(s) \| = 1 \quad \mbox{for} \quad \| W(0) \| =1 .
\end{equation}
\par
We derive a differential equation providing an approximation of
the learning system (\ref{output-signal}) with Oja's recurrence
rule (\ref{o-rule1}) in what follows. 
Equation (\ref{o-rule1}) takes the form
\begin{eqnarray}
\label{o-rule2}
W(s+1) =
\frac{W(s)+ \eta X(s) X^{T}(s) W(s)}{\| W(s) + \eta X(s) X^{T}(s) W(s)\|}
\end{eqnarray}
under the relation (\ref{output-signal}), which admits the Maclaurin
expansion
\begin{eqnarray}
\nonumber
W(s+1) &=& W(s) + \eta \{ X(s) X^{T}(s) W(s) \\
\label{Taylor}
& & - \left( W^{T}(s) X(s) X^{T}(s) W(s) \right) W(s) \} + O(\eta^{2})
\end{eqnarray}
if the leaning rate $\eta$ is sufficiently small. The $O(\eta^{2})$
in (\ref{Taylor}) denotes the second-order infinitesimal.
Elimination of the term $O(\eta^{2})$ from the rhs of (\ref{Taylor})
therefore provides us with the equation
\begin{eqnarray}
\nonumber
W(s+1) 
&=&
W(s) + \eta \big\{ X(s) X^{T}(s) W(s)
\\
\label{W}
&&
\qquad \qquad \quad - \left( W^{T}(s) X(s) X^{T}(s) W(s) \right) W(s) \big\} .
\end{eqnarray}
\subsection{Averaging}
We proceed an averaging of (\ref{W}) in what follows. Let $X(s)$ and
$W(s)$ be stochastic processes, which are statistically independent
to each other. On taking the expectation of (\ref{W}), we obtain
\begin{eqnarray}
\nonumber
&&
E[W(s+1) \, \vert \, W(s)] - E[W(s)] 
\\
\nonumber
&=&
\eta \big\{
E[X(s)X^{T}(s)] E[W(s)]
\\
\label{ex-o-rule}
&& \quad 
- \big(E[W(s)]^{T} E[X(s)X^{T}(s)] E[W(s)] \big) E[W(s)] \big\} ,
\end{eqnarray}
where the symbol $E[\cdot]$ denotes expectation operation.
On assuming the stochastic process $X(s)$ to be stationally,
the correlation matrix $E[X(s) X^{T}(s)]$ is kept invariant
along $s$, which is thereby diagonalized by an orthogonal
matrix $G$ to
\begin{eqnarray}
\label{C}
C = \mathrm{diag}(c_{1},c_{2},\cdots,c_{m}) = G^T E[X(s)X^{T}(s)]G.
\end{eqnarray}
Note that the orthogonal matrix $G$ do not depend on $s$,
so that all the eigenvalues $c_{j}$ of $E[X(s) X^{T}(s)]$
are kept invariant along $s$, too.
The change of variables
\begin{eqnarray}
\label{w}
w(t) = (w_{1}(t), w_{2}(t), \cdots, w_{m}(t))^{T}
= G^T W(t/\eta)
\end{eqnarray}
with the time-scaling
\begin{eqnarray}
\label{time-scaling}
t=\eta s
\end{eqnarray}
brings (\ref{ex-o-rule}) into the form,
\begin{equation}
\label{w-rec}
w(t + \eta) - w(t)
=
\eta \Big( Cw(t) - \big( w^{T}(t) C w(t) \big) w(t) \Big) .
\end{equation}
The differential equation
\begin{eqnarray}
\label{ALEH}
\frac{dw}{dt} 
= 
C w - (w^{T} C w ) w 
\end{eqnarray}
thereby emerges from (\ref{w-rec}) as a continuous-time approximation of
(\ref{o-rule2}) \cite{NHebb,Oja} in view of the stochastic approximation
theory \cite{Kush}. Throughout this paper, we will refer to (\ref{ALEH})
as the averaged learning equation of Hebb type (ALEH).
\subsection{The gradient-system form} 
We start with showing the norm preserving property
\begin{eqnarray}
\label{norm-preserve2}
\| w(t) \| =1 \quad \mbox{for} \quad \| w(0) \| =1
\end{eqnarray}
of the ALEH (\ref{ALEH}), which is understood to be a counterpart to
(\ref{norm-preserve}) of (\ref{output-signal}) with Oja's rule
(\ref{o-rule1}). Indeed, for any solution $w(t)$ of (\ref{ALEH}) with
$\| w(0) \| =1$, the calculation below shows (\ref{norm-preserve2});
\begin{eqnarray}
\label{wdw/dt}
\frac{d}{dt} \| w(t) \|^2
= 2w(t)^T \frac{dw}{dt}(t)
= 2(w(t)^{T} \,  C \, w(t)) (1-\| w(t) \|^2)=0.
\end{eqnarray}
Putting (\ref{w}), (\ref{time-scaling}) and (\ref{wdw/dt}) together,
we obtain
\begin{eqnarray}
\label{W-R}
\| W(s) \| = \| G w(\eta s) \| = \| w(\eta s) \| = \| w(0) \| = 1
\end{eqnarray}
for $\| W(0) \| = \| w(0) \| = 1$, 
as the counterpart to (\ref{norm-preserve}).
\par
Owing to the norm preserving property (\ref{norm-preserve2}),
we can restrict the ALEH on the $(m-1)$-dimensional unit sphere
\begin{eqnarray}
\label{Sph}
S^{m-1} = \Big\{ w \in {\bf R}^m \, \Big\vert \,
\| w \| = 1 \Big\} 
\end{eqnarray}
in ${\bf R}^{m}$. We endow $S^{m-1}$ with the
standard Riemannian metric
\begin{eqnarray}
\label{SPHmet}
\SPHmet{u}{u^{\prime}}_{w} = u^T u^{\prime}
\quad (u, u^{\prime} \in T_w S^{m-1}, \, w \in S^{m-1}),
\end{eqnarray}
where $T_w S^{m-1}$ denotes the tangent space of $S^{m-1}$ at
$w$ defined to be
\begin{eqnarray}
\label{tan-SPH}
T_{w} S^{m-1}
=\Big\{ u \in {\bf R}^m \, \Big\vert \, w^{T} u =0 \Big\}
\quad (w \in S^{m-1}). 
\end{eqnarray}
\par
According to Nakamura \cite{NHebb}, the ALEH on $S^{m-1}$ admits
the gradient-system form. Let the function,
\begin{eqnarray}
\label{Lambda}
\Lambda(w) = - \frac{1}{2} w^T C w 
= - \frac{1}{2}\sum_{k=1}^{m} c_{k} {w_{k}}^{2}
\quad (w \in S^{m-1}),
\end{eqnarray}
on $S^{m-1}$ be taken as the potential for the gradient-system form.
The gradient vector field $\mathrm{grad} \Lambda$
for the gradient system $(S^{m-1}, \SPHmet{\cdot}{\cdot}, \Lambda)$
is defined as follows: For a sufficiently small interval
$[a,b]$ with $a<0<b$, let us associate a smooth curve
$\gamma : [a,b] \rightarrow S^{m-1}$ with any $u \in T_w S^{m-1}$ in
the manner
\begin{eqnarray}
\label{gamma}
\tau \in [a, \, b] \mapsto \gamma (\tau ) \in S^{m-1} ,
\quad
\gamma (0)= w ,
\quad
\left. \frac{d\gamma}{d \tau } \right\vert_{\tau =0}= u \in T_{w}S^{m-1} . 
\end{eqnarray}
Then the gradient vector field $\mathrm{grad} \Lambda$ is defined
to satisfy \cite{KN}
\begin{eqnarray}
\label{def-grad-Lambda}
\SPHmet{\mathrm{grad} \Lambda (w)}{u}_w
=
\left. \frac{d}{d\tau } \right\vert_{\tau =0} \Lambda (\gamma (\tau ))
\quad (u \in T_{w}S^{m-1}).
\end{eqnarray}
Accordingly, $\mathrm{grad} \Lambda$ turns out to be
\begin{eqnarray}
\label{grad-Lambda}
\mathrm{grad} \, \Lambda  (w)
=
-C w + (w^{T} C w ) w ,
\end{eqnarray}
which conicides with the minus of the rhs of (\ref{ALEH}).
The ALEH on $S^{m-1}$ is thus written in the gradient-system form.
\section{Geometric devices}
This section provides geometric devices for realizing the ALEH
in the QIS.
\subsection{The QIS}
Following Uwano {\it et al}  \cite{UHI},
we introduce the quantum information space
(QIS), the space of regular density matrices endowed with the quantum SLD
(symmetric logarithmic derivertive) Fisher metric, in what follows.
\par
Let us consider the space of $m \times m$ regular density matrices
\begin{eqnarray}
\label{dot-P_m}
\Herpos = 
\{ \rho \in M(m,m) \, \vert \, 
\rho^{\dag} = \rho , \, \tr \rho =1 , \,
\rho : \mbox{positive definite} \} ,
\end{eqnarray}
where $M(m,m)$ denotes the set of $m \times m$ complex matrices.
The $\Herpos$ is endowed with the quantum SLD Fisher metric
$\QFmet{\cdot}{\cdot}$ as follows. 
\par
Let the tangent space of $\Herpos$ at $\rho$ be defined by
\begin{equation}
\label{tan-P_m} T_{\rho}\Herpos=\bigl\{\Xi\in M(m,m)\,\vert\,
\Xi^{\dagger}=\Xi,\,\tr\Xi=0\bigr\} \, .
\end{equation}
The symmetric logarithmic derivertive (SLD) to any tangent vector
$\Xi \in T_{\rho}\Herpos$ is defined to provide the Hermitean matrix
$\mathcal{L}_{\rho}(\Xi) \in M(m,m)$ subject to
\begin{equation}
\label{SLD} 
\frac{1}{2}\,\bigl\{\rho\mathcal{L}_{\rho}(\Xi)+
\mathcal{L}_{\rho}(\Xi)\,\rho\bigr\}=\Xi\qquad (\Xi\in
T_{\rho}\Herpos)\, .
\end{equation}
The quantum SLD Fisher metric, denoted by $\QFmet{\cdot}{\cdot}$,
is then defined to be
\begin{equation}
\label{defeq-Fisher}
\QFmet{\Xi}{\Xi^{\prime}}_{\rho}=\frac{1}{2}
\tr \left[\rho\bigl(L_{\rho}(\Xi)L_{\rho}(\Xi^{\prime})+
L_{\rho}(\Xi^{\prime})L_{\rho}(\Xi)\bigr)\right]
\qquad
(\Xi,\Xi^{\prime}\in T_{\rho}\Herpos) \, 
\end{equation}
(see also \cite{U,UHI,Brau,F,AN,H}). 
\par
We wish to present a more eplicit expression of $\QFmet{\cdot}{\cdot}$
in what follows. Let $\rho \in \Herpos$ be expressed as
\begin{eqnarray}
\begin{array}{l}
\rho=h\Theta h^{\dagger}, \quad h \in\mathrm{U}(m)
\\ \noalign{\medskip}
\label{rho}
\Theta=\mathrm{diag}(\theta_1,\ldots,\theta_m)
\quad
\mbox{with}
\quad
\tr \Theta =1 , \quad \theta_{k}>0 \,\, (k=1,2,\cdots m) ,
\end{array}
\end{eqnarray}
where $\mathrm{U}(m)$ denotes the group of $m \times m$ unitary
matrices. Expressing $\Xi \in T_{\rho} \Herpos$ as
\begin{eqnarray}
\label{Xi-chi}
\Xi = h \chi h^{\dagger}
\end{eqnarray}
with $h \in \mathrm{U}(m)$ in (\ref{rho}), we obtain an
explicit expression,
\begin{eqnarray}
\label{SLD-exp}
(h^{\dagger} {\cal L}_{\rho}(\Xi) h)_{jk}
= 
\frac{2}{\theta_j + \theta_k} \chi_{jk} \quad (j,k=1,2,\cdots m),
\end{eqnarray}
of the SLD to $\Xi \in T_{\rho}\Herpos$ \cite{UHI}.
Putting (\ref{rho})-(\ref{SLD-exp}) into (\ref{defeq-Fisher}), we have
\begin{equation}
\label{exp-QF}
\QFmet{\Xi}{\Xi^{\prime}}_{\rho}
=
2\sum_{j,k=1}^{m}
\frac{\overline{\chi}_{jk}\chi^{\prime}_{jk}}{\theta_j+\theta_k}
\end{equation}
where $\Xi^{\prime} \in T_{\rho}\Herpos$ is expressed as (\cite{UHI})
\begin{eqnarray}
\label{Xi'-chi'}
\Xi^{\prime} = h \chi^{\prime} h^{\dagger}.
\end{eqnarray}
\par
The space of $m \times m$ regular density matrices, $\Herpos$, endowed with
the quantum SLD Fisher metric $\QFmet{\cdot}{\cdot}$ defined above is
what we are referring to as the quantum information space (QIS) in the
present paper, which will be denoted also as the pair
$(\Herpos, \QFmet{\cdot}{\cdot})$ henceforth. 
\subsection{Metric preserving map}
We start with the Riemannian submanifold
\begin{eqnarray}
\label{cal-S}
{\cal S}_m 
=
S^{m-1} 
\backslash
\cup_{k=1}^{m} {\cal N}_m^{(k)} ,
\end{eqnarray}
with
\begin{eqnarray}
\label{N_m^k}
{\cal N}_m^{(k)}=
\{ w \in S^{m-1} \, \vert \, w_k=0 \}
\quad (k=1,2,\cdots ,m)
\end{eqnarray}
of $S^{m-1}$. From (\ref{cal-S}) and (\ref{N_m^k}), we immediately
obtain the coincidence
\begin{eqnarray}
\label{tan=}
T_{w}{\cal S}_m = T_{w} S^{m-1} \quad (w \in {\cal S}_m \subset S^{m-1})
\end{eqnarray}
of the tangent spaces $T_{w}{\cal S}_m$ and $T_{w} S^{m-1}$
(see also (\ref{tan-SPH})).
Then, the metric $\Smet{\cdot}{\cdot}$ of ${\cal S}_m$ is defined by
\begin{eqnarray}
\label{S-met}
\Smet{u}{u^{\prime}}_{w}
=
\SPHmet{\iota^{S}_{\ast ,w}(u)}{\iota^{S}_{\ast ,w}(u^{\prime})}_{w}
\quad
(u,u^{\prime} \in T_{w}{\cal S}_m),
\end{eqnarray}
where $\iota^{S}_{\ast , w}$ denotes the differential of the inclusion map
\begin{eqnarray}
\label{iota^S}
\iota^{S} : w \in {\cal S}_m \mapsto w \in S^{m-1}
\end{eqnarray}
at $w \in {\cal S}_m$ (see Appendix~A for the differential of maps).
According to (\ref{iota^S}) and Appendix~A, $\iota^S_{\ast ,w}(u)$ turns
out to be
\begin{eqnarray}
\label{iota^S_*}
\iota^S_{\ast ,w} (u) = u \in T_{w} \Herpos
\quad (u \in T_{w}{\cal S}_m) .
\end{eqnarray}
\par
We move to consider the Riemannian submanifold
\begin{eqnarray}
\nonumber
&&
{\cal D}_m
=
\Big\{ \Theta \in \Herpos \,\, \Big\vert \,\,
\Theta = \mathrm{diag} \, (\theta_{1}, \cdots , \theta_{m}) , 
\\
\label{cal-D}
&&
\phantom{{\cal D}_m=\Big\{ \rho \in \Herpos \,\, \Big\vert \,\,}
\sum_{k=1}^{m}\theta_k =1, \,\, \theta_{k} >0 \,\,\, (k=1,2,\cdots ,m)
\Big\} 
\end{eqnarray}
of the QIS $(\Herpos, \QFmet{\cdot}{\cdot})$. The tangent space
$T_{\Theta}{\cal D}_m$ at $\Theta$ takes the form
\begin{eqnarray}
\label{tan-cal-D}
T_{\Theta}{\cal D}_m
=
\Big\{ Z \in M(m,m) \, \Big\vert \, 
Z =  \mathrm{diag} \, (\zeta_{1}, \cdots , \zeta_{m}),
\,
\sum_{j=1}^{m} \zeta_j = 0 \Big\} \, \subset T_{\Theta}\Herpos .
\end{eqnarray}
The metric $\Dmet{\cdot}{\cdot}$ of ${\cal D}_m$ is defined by
\begin{eqnarray}
\label{QF-D}
\Dmet{Z}{Z^{\prime}}_{\Theta}
=
\QFmet
{\iota^D_{\ast , \Theta}(Z)}
{\iota^D_{\ast , \Theta}(Z^{\prime})}_{\Theta} \quad
(Z,Z^{\prime} \in T_{\Theta}{\cal D}_m),
\end{eqnarray}
where $\iota^D_{\ast , \Theta}$ denotes the differential
of the inclusion map
\begin{eqnarray}
\label{iota^D}
\iota^D : \Theta \in {\cal D}_m \mapsto \Theta \in  \Herpos 
\end{eqnarray}
 at $\Theta \in{\cal D}_m$ (see Appendix~A).
Like in the case of $\iota^S_{\ast ,w}$ (see (\ref{iota^S_*})),
$\iota^D_{\ast ,\Theta} (Z)$ turns out to be
\begin{eqnarray}
\label{iota^D_*}
\iota^D_{\ast ,\Theta} (Z) = Z \in T_{\Theta} \Herpos
\quad (Z \in T_{\Theta}{\cal D}) ,
\end{eqnarray}
\par
Under the preparation above, we consider the smooth map $\mu$
of ${\cal S}_m$ to ${\cal D}_m$ of the form  
\begin{eqnarray}
\label{mu}
\mu : w = (w_1, w_2,\cdots,w_m)^T \in {\cal S}_m
\mapsto
\mathrm{diag} \, (w_1^2, w_2^2, \cdots, w_m^2) \in {\cal D}_m ,
\end{eqnarray}
whose differential of the map $\mu_{\ast ,w}$ at $w$ takes the form
\begin{eqnarray}
\label{mu_*-exp}
\mu_{\ast ,w}(u)
=
2 \, \mathrm{diag} \,  (w_1u_1 , w_2u_2, \cdots , w_mu_m)
\quad (u \in T_w{\cal S}_m).
\end{eqnarray}
We have the following lemma for $\mu$.
\begin{lemma}\label{lemma-mu}
The map $\mu$ satifies (i)-(v) in the following:
\par\noindent
(i) The $\mu$ is surjective.
\par\noindent
(ii) The $\mu$ is metric-preserving up to the constant multiple $4$;
\begin{eqnarray}
\label{met-pres}
4 \, \Smet{u}{u^{\prime}}_w
=
\Dmet{\mu_{\ast ,w}(u)}{\mu_{\ast ,w}(u^{\prime})}_{\mu (w)}
\quad 
(w \in {\cal S}_m , \, u, u^{\prime} \in T_w {\cal S}_m)
\end{eqnarray}
holds true, where $T_w {\cal S}_m$ is defined by (\ref{tan=}) with
(\ref{tan-SPH}).
\par\noindent
(iii) Under the $({\bf Z}_2)^m$-action on ${\cal S}_m$ defined by
\begin{eqnarray}
\nonumber
&&
\phi_{\sigma} : w=(w_1,w_2, \cdots, w_m)^T \in {\cal S}_m
\\
\label{phi}
&&
\qquad \mapsto
(\sigma_1 w_1, \sigma_2 w_2, \cdots , \sigma_m w_m)^T \in {\cal S}_m
\end{eqnarray}
with
\begin{eqnarray}
\label{Z_2^m}
\sigma =(\sigma_1, \sigma_2, \cdots , \sigma_m)^T  \in ({\bf Z}_2)^m
\quad ({\bf Z}_2 =\{ -1,1 \} ),
\end{eqnarray}
the $\mu$ is invariant;
\begin{eqnarray}
\label{mu-inv}
\mu ( \phi_{\sigma} (w))= \mu (w) \quad
(w \in {\cal S}_m , \, \sigma \in ({\bf Z}_2)^m).
\end{eqnarray}
(iv) The coincidence $\mu (w) = \mu (w^{\prime})$ holds true if and only if
there exists a certain $\sigma \in ({\bf Z}_2)^m$ subject to
$w^{\prime}=\phi_{\sigma}(w)$.
\par\noindent
(v) The restrict of $\mu$ to each of
\begin{eqnarray}
\label{pos-part}
{\cal S}_m^{\sigma} = \{ w \in {\cal S}_m \, \vert \,
\sigma_j w_j >0 \, (j=1,2,\cdots,m) \} \quad (\sigma \in ({\bf Z}_2)^m)
\end{eqnarray}
is diffeomorophic (smooth, injective and surjective) of ${\cal S}_m^{\sigma}$.
\end{lemma}
The items other than (ii) in Lemma \ref{lemma-mu} are proved easily by
straightforward calculations. The proof of (ii) is consigned in Appendix~B.
\begin{remark}\label{rem-cal-S}
If we apply (\ref{mu}) formally to any $w \in {\cal N}_m^{(k)}$,
$\mu (w)$ is out of $\Herpos$. This is the account for eliminating
${\cal N}_m^{(k)}$s from $S^{m-1}$ to consider ${\cal S}_m$.
\end{remark}
\begin{remark}\label{imbedding}
From (\ref{mu_*-exp}) and (\ref{tan=}) with (\ref{tan-SPH}),
we have $\mathrm{rank}\, \mu_{\ast ,w}=m-1 = \mathrm{dim} \,{\cal D}_m$.
The map $\mu$ is hence called an immersion \cite{KN}.
\end{remark} 
\begin{remark}\label{quotient}
In view of (i) and (iv) in Lemma \ref{lemma-mu}, ${\cal D}_m$ is understood
as the quotient space ${\cal S}_m / ({\bf Z}_2)^m$ of ${\cal S}_m$.
Further, combining (v) with (i) and (iv), we see that ${\cal S}_m^{\sigma}$
with every $\sigma \in ({\bf Z}_2)^m$ is diffeomorphic to ${\cal D}_m$.
\end{remark}
\section{The ALEH on the QIS}
We are now in a position to construct the gradient system on the QIS
realizing an extention of the ALEH by making full use uses of the
geometric devices developped in section~3.
\subsection{Mapping the ALEH through $\mu_{\ast ,w}$}
In this subsection, the vector field for the ALEH (\ref{ALEH}) on
${\cal S}_m$ is shown to be mapped to a vector field on ${\cal D}_m$.
Note that the fact above is not so trivial since the map $\mu$ is not
injective (see (ii) of Lemma \ref{lemma-mu} and Remark \ref{quotient}).
\par
In view of (i) of Lemma~\ref{lemma-mu} and Remark~\ref{rem-cal-S},
we start with the restrict of the ALEH (\ref{ALEH}) to ${\cal S}_m$.
Since $\mathrm{grad} \Lambda $ satisfies
\begin{eqnarray}
\label{grad-bd}
( \mathrm{grad} \Lambda )(w) \in T_w {\cal N}_m^{(k)}
\quad \mbox{for} \quad w \in {\cal N}_m^{(k)}
\quad (k=1,2,\cdots , m) ,
\end{eqnarray}
all trajectories with initial condition
$w(0) \in {\cal N}_m^{(k)}$ are confined in ${\cal N}_m^{(k)}$
($k=1,2,\cdots ,m$), respectively. On the same account, all the trajectories
with $w(0) \notin \cup_{k=1}^m {\cal N}_m^{(k)}$
never intersect with $\cup_{k=1}^m {\cal N}_m^{(k)}$.
Accordingly, the ALEH (\ref{ALEH}) can be dealt with
on ${\cal S}_m = S^{m-1} \backslash \cup_{k=1}^m {\cal N}_m^{(k)}$.
\par
The next thing to show is the $({\bf Z}_2)^m$-invariance of
$\mathrm{grad} \Lambda$, which is equivalent for
\begin{eqnarray}
\label{grad-equiv}
\mathrm{grad} \Lambda (\phi_{\sigma} (w))
=
(\phi_{\sigma})_{\ast ,w} ( \mathrm{grad} \Lambda (w))
\quad
(w \in {\cal S}_m, \, \sigma \in ({\bf Z}_2)^m)
\end{eqnarray}
(see (\ref{phi}) and (\ref{Z_2^m}) for $\phi_{\sigma}$).
Indeed, on taking the expression
\begin{eqnarray}
\label{phi_*}
(\phi_{\sigma})_{\ast ,w}(u)
=
(\sigma_1 u_1, \sigma_2 u_2, \cdots , \sigma_m u_m)^T
\quad (u \in T_w {\cal S}, \, \sigma \in ({\bf Z}_2)^m)
\end{eqnarray}
into account (see Appendix~A for the differential of maps),
Equation (\ref{grad-equiv}) is shown to hold true by a simple calculation.
\par
The $({\bf Z}_2)^m$-invariance thus shown is put together with (iv) of
Lemma \ref{lemma-mu} to ensure the existence of the vector
field $\mu_{\ast} \mathrm{grad} \Lambda$ on ${\cal D}_m$
subject to
\begin{eqnarray}
\label{mu-grad-Lambda}
\mu_{\ast} \mathrm{grad} \Lambda (\mu (w))
=
\mu_{\ast ,w} (\mathrm{grad} \Lambda (w))
\quad (w \in {\cal S}_m).
\end{eqnarray}
Note that $\mu_{\ast} \mathrm{grad} \Lambda $ is well-defined since
the equation
\begin{eqnarray}
\nonumber
&&
\mu_{\ast , \phi_{\sigma}w} ( \mathrm{grad} \Lambda ( \phi_{\sigma}w))
=
\mu_{\ast , \phi_{\sigma}w} ( (\phi_{\sigma})_{\ast , w}
(\mathrm{grad} \Lambda (w))
\\
\label{well-def}
&&
\phantom{
\mu_{\ast , \phi_{\sigma}w} ( \mathrm{grad} \Lambda ( \phi_{\sigma}w))
}
=\mu_{\ast , w} ( \mathrm{grad} \Lambda (w))
\quad (w \in {\cal S}_m, \, \sigma \in ({\bf Z}_2)^m)
\end{eqnarray}
holds true. The first equality follows from (\ref{grad-equiv})
and the second one from
\begin{eqnarray}
\label{compo-diff}
\mu_{\ast , \phi_{\sigma}w} \circ (\phi_{\sigma})_{\ast , w}
=
(\mu \circ \phi_{\sigma})_{\ast , w}
=
\mu_{\ast ,w}
\end{eqnarray}
with (\ref{mu-inv}). Equation (\ref{compo-diff}) follows from
(\ref{mu_*-exp}) with (\ref{phi_*}) immediately.
\begin{remark}\label{interpretation}
On closing this subsection, we give a naive description of
the vector field $\mu_{\ast} \mathrm{grad} \Lambda$ on ${\cal D}_m$.
Recalling Remark \ref{quotient}, we may understand that the image
${\cal D}_m$ of the map $\mu$ is a $2^m$-fold copy of ${\cal S}_m$.
Therefore, in view of (\ref{well-def}), the vector field
$\mu_{\ast} \mathrm{grad} \Lambda$ on ${\cal D}_m$ can be regarded as a 
$2^m$-folded copy of the ALEH on ${\cal S}_m$.
\end{remark}
\subsection{The gradient system on the QIS realizing the ALEH}
We are now in a position to seek a gradient system on the QIS realizing
the ALEH. Namely, what we we are to seek is a gradient system
that realizes
\begin{eqnarray}
\label{grad=}
\mathrm{grad} L (\Theta) = \mu_{\ast} \mathrm{grad} \Lambda (\Theta)
\quad (\Theta \in {\cal D}_m) ,
\end{eqnarray}
where $L$ denotes the potential for the gradient system.
We note here that the gradient vector field $\mathrm{grad} L$
associated with the potential $L$ is defined on the QIS
$(\Herpos , \QFmet{\cdot}{\cdot})$ in the following way
(cf. (\ref{gamma}) and (\ref{def-grad-Lambda})).
For a sufficiently small interval
$[a,b]$ with $a<0<b$, let us associate a smooth curve
$r : [a,b] \rightarrow \Herpos$ with any $\Xi \in T_{\rho} \Herpos$ in
the manner
\begin{eqnarray}
\label{r}
\tau \in [a, \, b] \mapsto r (\tau ) \in \Herpos ,
\quad
r (0)= \rho ,
\quad
\left. \frac{dr}{d \tau } \right\vert_{\tau =0}= \Xi \in T_{\rho} \Herpos . 
\end{eqnarray}
Then the gradient vector field $\mathrm{grad} L$ is defined
to satisfy \cite{KN}
\begin{eqnarray}
\label{def-grad-L}
\QFmet{\mathrm{grad} L ( \rho )}{\Xi}_{\rho}
=
\left. \frac{d}{d\tau } \right\vert_{\tau =0} L (r (\tau ))
\quad (\Xi \in T_{\rho}\Herpos ).
\end{eqnarray}
Instead of (\ref{grad=}), we are to consider a weaker condition than
(\ref{grad=}) below in order that we can fix {\it a candidate} easily
for the potential. 
The weaker condition to be dealt with is
\begin{eqnarray}
\label{part=}
\QFmet{\mathrm{grad} L (\Theta)}{\iota^D_{\ast ,\Theta}Z}_{\Theta}
=
\Dmet{\mu_{\ast} \mathrm{grad} \Lambda (\Theta)}{Z}_{\Theta}
\quad
(\Theta \in {\cal D}_m , Z \in T_{\Theta} {\cal D}_m) .
\end{eqnarray}
We show the following Lemma.
\begin{lemma}\label{lemma-L}
If the potential $L$ for a gradient system
$(\Herpos , \QFmet{\cdot}{\cdot}, L)$ satisfies
\begin{eqnarray}
\label{L-Lambda}
L(\mu (w))= 4 \Lambda (w) \quad (w \in {\cal S}_m), 
\end{eqnarray}
then (\ref{part=}) holds true.
\end{lemma}
{\it Proof}\quad
It follows from (\ref{mu}) that the unique $w \in {\cal S}_m^{\sigma}$
with $\sigma = (1,1,\cdots ,1)^T$ subject to $\mu (w) = \Theta$,
where ${\cal S}_m^{\sigma}$ is defined by (\ref{pos-part}).
Further, recalling (v), we can take the unique $u \in T_{w}{\cal S}_m$
subject to $\mu_{\ast ,w}u=Z$. Then for a sufficiently small interval
$[a,b]$ with $a<0<b$, we consider a curve $\gamma (\tau)$
subject to (\ref{gamma}) and a curve $r(t)=\mu (\gamma (\tau))$.
The setting above is put together with (\ref{def-grad-Lambda}), 
(\ref{met-pres}), (\ref{mu-grad-Lambda}), (\ref{def-grad-L}) and
the assumption (\ref{L-Lambda}) to show
\begin{eqnarray}
\nonumber
& \phantom{\! = \!} &
\QFmet{\mathrm{grad} L (\Theta)}{\iota^D_{\ast ,\Theta}Z}_{\Theta}
\\
\nonumber
& \! = \! &
\QFmet{\mathrm{grad} L (\Theta)}{Z}_{\Theta}
=
\left. \frac{d}{d\tau} \right\vert_{\tau =0} L(r(\tau ))
=
\left. \frac{d}{d\tau} \right\vert_{\tau =0} L(\mu (\gamma (\tau )))
\\
\nonumber
& \! = \! &
\left. \frac{d}{d\tau} \right\vert_{\tau =0} 4\Lambda (\gamma (\tau ))
=
4 \Smet{ \mathrm{grad} \Lambda (w)}{u}_{w}
\\
\label{lhs=rhs}
& \! = \! &
\Dmet{\mu_{\ast ,w}(\mathrm{grad} \Lambda (w))}
{\mu_{\ast ,w}(u)}_{\Theta}
=
\Dmet{\mu_{\ast} \mathrm{grad} \Lambda (\Theta)}{Z}_{\Theta}.
\end{eqnarray}
This completes the proof.
\par\smallskip
Owing to Lemma \ref{lemma-L}, we can choose
\begin{eqnarray}
\label{L}
L(\rho) = -2 \, \mathrm{tr} \, ( C \rho )
\end{eqnarray}
as a candidate for the potential realizing (\ref{grad=}),
where $C$ is the diagonal matrix given in (\ref{C}).
Note that we can confirm (\ref{L-Lambda}) for $L$ of (\ref{L})
by the calculation,
\begin{eqnarray}
\label{confirm}
L(\mu (w))= - 2 \, \mathrm{tr} \, ( C \mu (w) )
= -2 \sum_{j=1}^{m} c_j w_j^2 = 4 \Lambda (w),
\end{eqnarray}
made with (\ref{Lambda}) and (\ref{mu}). We move to draw the gradient
equation for $L$ of (\ref{L}) along with the framework of gradient systems
on the QIS developped by the authors \cite{UY,UY2}. Following to
\cite{UY,UY2}, we define the Hermitean matrix ${\cal M}(L)$ for $L$
to be
\begin{eqnarray}
\label{cal-M}
\big( {\cal M}(L) \big)_{jk}
=
\left\{ \begin{array}{ll}
\displaystyle{
\frac{\partial L}{\partial \overline{\rho}_{jk}}
=
\overline{ \Big( \frac{\partial L}{\partial \rho_{jk}} \Big) }
}
& \quad (1 \leq j < k \leq m)
\\
\noalign{\smallskip}
\displaystyle{
\frac{\partial L}{\partial \rho_{jj}}
}
& \quad (j=k=1,2, \cdots ,m) ,
\end{array}
\right.
\end{eqnarray}
where the partial differentiations, $\partial /\partial \rho_{jk}$
and $\partial /\partial \overline{\rho}_{jk}$, stand for
\begin{eqnarray}
\label{diff-rho}
\begin{array}{l}
\displaystyle{
\frac{\partial}{\partial \rho_{jk}}
=
\frac{1}{2}\Big( 
\frac{\partial}{\partial \xi_{jk}} -i \frac{\partial}{\partial \eta_{jk}}
\Big) 
}
\\ \noalign{\medskip}
\displaystyle{
\frac{\partial}{\partial {\overline \rho}_{jk}}
=
\frac{1}{2}\Big( 
\frac{\partial}{\partial \xi_{jk}} + i \frac{\partial}{\partial \eta_{jk}}
\Big)
}
\end{array}
\quad (1 \leq j<k \leq m)
\end{eqnarray}
with
\begin{eqnarray}
\label{xi-eta}
\xi_{jk} = \Re (\rho_{jk}) , \quad
\eta_{jk} = \Im (\rho_{jk}) \quad
\quad (1 \leq j<k \leq m).
\end{eqnarray}
The symbols $\Re$ and $\Im$ indicate the real part and the imaginary one,
respectively.
In contrast, the $\rho_{jj}$s are thought of as real variables to give
rise to $\partial /\partial \rho_{jj}$s in usual way. 
In terms of ${\cal M}(L)$, the gradient equation for $L$ is written as
\begin{eqnarray}
\label{grad-L-M}
\mathrm{grad} L ( \rho )
=
\frac{1}{2} \Big( 
\rho {\cal M}(L) + {\cal M}(L) \rho \Big) 
- \Big( \tr \big(\rho {\cal M}(L) \big) \Big) \rho
\end{eqnarray}
(see \cite{UY,UY2}). Since we have
\begin{eqnarray}
\label{ML}
{\cal M}(L) = -2C
\end{eqnarray}
through a straightforward calculation of (\ref{cal-M}) with (\ref{L}),
the gradient vector field for $L$ is exressed as
\begin{eqnarray}
\label{grad-L-fin}
\mathrm{grad} L ( \rho )
=
-(\rho C + C \rho ) + 2 \Big( \tr \big(\rho C\big) \Big) \rho 
\end{eqnarray}
with the diagonal matrix $C$ of (\ref{C}). To summarize, we have the
following.
\begin{lemma}\label{grad-sys-L}
The gradient system $(\Herpos , \QFmet{\cdot}{\cdot}, L)$ is governed
by the differential equation of motion
\begin{eqnarray}
\label{grad-eq-L}
\frac{d\rho}{dt}
=
(\rho C + C \rho ) - 2 \Big( \tr \big(\rho C\big) \Big) \rho ,
\end{eqnarray}
where $C$ is the diagonal matrix of (\ref{C}).
\end{lemma}
We are at the final stage to show (\ref{grad=}) for $L$ of (\ref{L}).
On combining (\ref{grad-Lambda}) with (\ref{mu_*-exp}), 
the rhs of (\ref{grad=}) is caluculated to be
\begin{eqnarray}
\nonumber
& \phantom{\! = \!} &
\mu_{\ast ,w} ( \mathrm{grad} \Lambda (w)) 
\\ \noalign{\medskip}
\label{mu_*-grad-Lambda-fin}
&\! = \!&
2 \, 
\mathrm{diag} \, (
w_1 (\mathrm{grad} \Lambda (w))_1,
w_2 (\mathrm{grad} \Lambda (w))_2, \cdots , 
w_m (\mathrm{grad} \Lambda (w))_m )
\end{eqnarray}
with
\begin{eqnarray}
\label{grad-Lambda_j}
(\mathrm{grad} \Lambda (w))_j
=
-c_j w_j + \Big( \sum_{k=1}^m c_k w_k^2 \Big) w_j
\quad (j=1,2,\cdots ,m),
\end{eqnarray}
while (\ref{grad-L-fin}) is put together with (\ref{mu}) to show
\begin{eqnarray}
\nonumber
\mathrm{grad} L ( \mu (w) )
&\! = \! &
- 2 \, \mathrm{diag} (c_1w_1^2, c_2w_2^2, \cdots , c_mw_m^2)
\\
\label{grad-L-D}
&&
+ 2 \Big( \sum_{k=1}^m c_k w_j^2 \Big) \,
\mathrm{diag} (c_1w_1^2, c_2w_2^2, \cdots , c_mw_m^2)
\end{eqnarray}
for $w \in {\cal S}_m$. 
Equations (\ref{mu_*-grad-Lambda-fin})-(\ref{grad-L-D}) therefore
confirms (\ref{grad=}) for $L$ of (\ref{L}).
Finally, we reach to the following theorem.
\begin{theorem}
\label{ALEH-QIS}
The gradient system $(\Herpos , \QFmet{\cdot}{\cdot}, L)$ governed
by the equation of motion (\ref{grad-eq-L}) realizes the ALEH in the
QIS in the following sense:
In the manner of (\ref{grad=}),
the gradient vector field $\mathrm{grad} L$ realizes the vector field
mapped through (\ref{mu-grad-Lambda}) from the gradient vector field
describing the ALEH on ${\cal S}_m$.
\end{theorem}
\section{Concluding remarks}
We have constructed successfully the gradient system
on the QIS (GS-QIS) which realizes the ALEH on the
submanifold ${\cal D}_m$. The success is due to the geometric
devices developped in Section 3: Especially, a key fact is
that the immersion $\mu$ of ${\cal S}_m$ to ${\cal D}_m$
is put together with the $({\bf Z}_2)^m$-symmetry of the ALEH
to realize the $2^m$-folded copy of the ALEH on ${\cal D}_m$.
\par
The behavior of trajectories of the GS-QIS not on ${\cal D}$ is open still,
though the trajectories on ${\cal D}$ are understood by Nakamura \cite{NHebb}:
The GS-QIS is expected to have a \lq global convergence' property.
\par
Integrability of the GS-QIS is an open question, too.
However, in the special case of $C=I$, integrability is allowed
like in the cases of \cite{U} and \cite{UY,UY2} with $C=-2I$,
since the $\mathrm{U}(m)$ action,
\begin{eqnarray}
\label{U(m)}
\rho \in \Herpos \mapsto h \rho h^{\dagger} \in \Herpos
\quad (h \in \mathrm{U}(m)),
\end{eqnarray}
is allowed to be a symmetry, whose role is studied in \cite{U}. 
\par
Relation of the gradient systems obtained in \cite{U,UHI,UY,UY2}
to physical systems would be a big problem. The horizontal lift of
gradient vectors to those on the space of ordered-tuples of multi-qubit
states developped in \cite{U,UHI} could be a clue, and the work of
Braunstein \cite{Brau} another one. Further, the authors might say that
the complementary property between gradient vector fields and the Hamiltonian
ones in classical mechanics is worth taken together with the horizontal lift.  
\par\bigskip\noindent
{\bf Acknowlegement}
\quad
The authors thank Professor Toshihiro Iwai at Kyoto University
for his comments on the present work.
This work is partly supported by Special Research Funds,
A2 (2008) and B14 (2009), Future University Hakodate.
\appendix
\setcounter{section}{0}
\section{Differential of maps} 
With the aim of a brief and clear description, we restric to our attention
to smooth maps among spaces in $M(m,m)$, the space of $m \times m$ complex
matrices. Let $M_1$ and $M_2$ be manifolds in $M (m,m)$ and $\psi$ a map
of $M_1$ to $M_2$;
\begin{eqnarray}\label{psi}
\psi : M_1 \rightarrow M_2.
\end{eqnarray}
For a given point $p \in M_1$, the differential $\psi_{\ast ,p}$ of $\psi$
at $p$ is defined as follows. Like (\ref{gamma}), for a sufficiently small
interval $[a,b]$ $a<0<b$, let us associate a sommth curve
$q : [a,b] \rightarrow M_1$
\begin{eqnarray}
\label{zeta}
\tau \in [a, \, b] \mapsto q (\tau ) \in M_1 ,
\quad
q (0)= w ,
\quad
\left. \frac{dq}{d \tau } \right\vert_{\tau =0}= w \in T_{p}M_1 . 
\end{eqnarray}
On using the smooth curve $q (\tau)$ introduced above, 
the differential $\psi_{\ast ,p}$ of the map $\psi$ at $p$ is defined to be
\begin{eqnarray}
\label{psi_*}
\psi_{\ast ,p}(v) 
=
\left. \frac{d}{d \tau} \right\vert_{\tau =0} \psi (q (\tau))
\quad (v \in T_{p}M_1).
\end{eqnarray}
(see \cite{KN}, for example).
The differential maps, $\iota^S{\ast, w}$, $\iota^D_{\ast , \Theta}$,
$\mu_{\ast ,w}$ and $(\phi_{\sigma})_{\ast ,w}$ are defined by
(\ref{zeta}) and (\ref{psi_*}) with $({\cal S}_m, S^{m-1} , w , \iota^S )$,
$({\cal D}_m, \Herpos , \Theta , \iota^D )$,
$({\cal S}_m, {\cal D}_m , w, \mu )$ and
$({\cal S}_m, {\cal S}_m , w, \psi_{\sigma} )$
in place of $(M_1 , M_2, p, \psi )$.
\section{Proof of Lemma \ref{lemma-mu}}
We show (ii) of Lemma \ref{lemma-mu}. On setting
$h=id$, $\Xi=Z$ and $\Xi^{\prime}=Z^{\prime}$ in
(\ref{Xi-chi}) and (\ref{Xi'-chi'}), Equation (\ref{exp-QF})
is put together with (\ref{QF-D}) to show
\begin{eqnarray}
\label{rhs-lemma-mu}
\Dmet{Z}{Z^{\prime}}_{\Theta}
=
\sum_{k=1}^m \frac{\zeta_k \zeta^{\prime}_k}{\theta_k}.
\end{eqnarray}
Further the substitutions $Z=\mu_{\ast ,w}(u)$ and
$Z^{\prime}=\mu_{\ast ,w}(u^{\prime})$ with (\ref{mu_*-exp}) bring us
to have
\begin{eqnarray}
\label{rhs-calc}
\Dmet{Z}{Z^{\prime}}_{\Theta}
=
\sum_{k=1}^m \frac{\zeta_k \zeta^{\prime}_k}{\theta_k}.
=
\sum_{k=1}^m \frac{(2w_k u_k) (2w_k u^{\prime}_k)}{w_k^2}
=
4 \sum_{k=1}^m u_k u^{\prime}_k
=
4 \Smet{u}{u^{\prime}}_{w}.
\end{eqnarray}
This ends the proof.

\end{document}